\documentclass{article}
\usepackage{graphicx} 
\usepackage[english]{babel}
\usepackage{amsmath, amsfonts, amssymb, amsthm}
\usepackage{parskip}
\usepackage{mathtools}
\usepackage{tikz}
\usetikzlibrary{positioning}
\usetikzlibrary{calc}
\usepackage{caption}
\usepackage{authblk}
\usepackage[super]{nth}

\usepackage{comment}

\usepackage[labelfont=bf]{caption}
\newtheorem{theorem}{Theorem}[subsubsection]
\newcounter{theoremaux}
\newtheorem{proposition}{Proposition}[subsection]
\newtheorem{corollary}{Corollary}[subsection]
\newtheorem{lemma}{Lemma}[subsubsection]
\newcounter{lemmaaux}
\theoremstyle{definition}
\newtheorem{definition}{Definition}[subsubsection]
\newcounter{definitionaux}
\newtheorem{conjecture}{Conjecture}[section]

\title{Babai Numbers and Babai Spectra of Paths and Cycles}
\author[(1)]{Peter Johnson}
\author[(2)]{Celalettin Kaya}
\author[(3)]{Ryan W. Matzke}
\affil[(1)]{\normalsize{Department of Mathematics and Statistics, Auburn University, Auburn, \makebox{AL 36849,} USA, johnspd@auburn.edu}}
\affil[(2)]{\normalsize{Department of Mathematics, Çankırı Karatekin University, Çankırı, \linebreak 18100, Turkey, ckaya@karatekin.edu.tr}}
\affil[(3)]{\normalsize{Department of Mathematics, Vanderbilt University, Nashville, \newline TN 37240, USA, ryan.w.matzke@vanderbilt.edu \linebreak Research by this co-author was funded by the NSF Research in \linebreak Mathematical Sciences Postdoctoral Fellowship 2202877}}
\date{August 2024}

\begin{document}

\maketitle

\begin{abstract}
\noindent We study Babai numbers and Babai $k$-spectra of paths and cycles. We completely determine the Babai numbers of paths $P_n$ for $n>1$ and $1 \leq k \leq n-1$, and the Babai $k$-spectra for $P_n$ when $1 \leq k \leq n/2$ . We also completely determine Babai numbers and Babai $k$-spectra of all cycles $C_n$ for $k \in \{1,2\}$ and $n \geq 3$ if $k=1$ and $n > 3$ if $k=2$.   
\end{abstract}


\renewcommand{\thefootnote}{\fnsymbol{footnote}}

\section{Introduction}

\counterwithin{lemma}{section}
\counterwithin{theorem}{section}
\counterwithin{definition}{section}

\hspace*{10pt}The Babai numbers were proposed by Laszlo Babai in conversation with the first author of \cite{aj}, where the numbers may have first appeared in print. (We hear that Paul Erd\H{o}s had the same idea, but in a limited context.) Before giving the definition of Babai numbers, we first need some basic notions, also to be found in \cite{aj} and \cite{jk}.

A \textit{distance function} on a non-empty set $X$ is a function $\rho : X \times X \to [0, \infty)$ such that, for all $x,y \in X$, (i) $\rho(x,y) = \rho(y,x)$ and (ii) $\rho(x,y) = 0 \Leftrightarrow x = y$. If, in addition, (iii) $\rho(x,y) \leq \rho(x,z)+\rho(z,y)$ for all $x,y,z \in X$, then $\rho$ is a \textit{metric}, and the pair $(X,\rho)$ is a \textit{metric space}.

A \textit{coloring} of a set $X$ is a function $\varphi: X \to C$ for some set $C$, called the set of colors; when $C$ is to be finite, we usually take $C$ to be a set of positive integers, that is, $C=[m]=\{1,2,\ldots, m\}$ for some $m \in \mathbb{Z}^+$. Such a coloring is completely described by the color sets $\varphi^{-1}(\{c\}) = \{x \in X \mid \varphi(x) = c\}$ for $c \in C$. Some of these sets may be empty. If we allow empty sets in a partition, then the color sets partition $X$. Conversely, every partition of $X$ constitutes a coloring of $X$.

Suppose that $\rho$ is a distance function on a set $X$, $\varphi: X \to C$ is a coloring of $X$, $d > 0$ and $c \in C$. We will say that the distance $d$ is \textit{forbidden for the color $c$} if $\varphi^{-1}(\{c\})$ contains no two points $x, y \in X$ such that $\rho(x, y) = d$. To put it another way, if $x, y \in X$ and $\rho(x, y) = d$, then $x$ and $y$ cannot both be colored $c$. If the distance $d$ is forbidden for every $c \in C$, then we say that the distance $d$ is \textit{forbidden by the coloring $\varphi$}.

Suppose that $X$ and $\rho$ are as above, and $D \subseteq (0, \infty)$. The \textit{distance graph} $G(X, D)$ is the simple graph with vertex set $X$ and in which $x, y \in X$ are adjacent if and only if $\rho(x, y) \in D$. (It is conventional to suppress mention of $\rho$ in the notation. In this paper, it will always be clear what distance function is under discussion.) When $D = \{d\}$, we can use the notation $G(X, d)$ as well as $G(X, \{d\})$.

If $P$ is a graph parameter, it is customary to shorten $P(G(X, D))$ to $P(X, D)$. Thus, the chromatic number of $G(X, D)$ is denoted $\chi(X, D)$. Due to the way adjacency is defined in $G(X, D)$, we have $\chi(X, D) = min[\vert C \vert :$ there is a coloring $\varphi : X \to C$ which forbids each $d \in D]$.

For a given $(X, \rho)$, let $R(X, \rho)=\rho(X \times X) \setminus \{0\}$. That is, $R(X, \rho)$ is the set of positive distances realized by $\rho$, between points of $X$.

Now, we can give the definition of Babai numbers.

\begin{definition}
Let $k$ be a positive integer satisfying $k \leq \vert R(X, \rho) \vert$. The \textit{$k^{th}$ Babai number of $(X, \rho)$} is $B_k(X) = sup[\chi(X, D) : D \subseteq R(X, \rho) \text{ and } \vert D \vert = k]$.
\end{definition}

Note that since $\chi(X, D) \leq \vert X \vert$ for all $D \subseteq (0, \infty)$, the sup, or least upper bound, in the definition of $B_k(X)$ is valid. However, in most cases, including all of those of interest to us in this paper, $B_k(X)$ is finite and the sup above is a maximum.

The definition of $B_k(X)$ given here is equivalent to the original definition given in \cite{aj} (except that the erroneous ``max'' there should have been ``sup''), but there the sup is taken over all $k$-subsets of $(0, \infty)$, whereas here the sup is taken over all $k$-subsets of $R(X, \rho)$. We have introduced this complication in preparation for the introduction of Babai $k$-spectra, to be defined below. 

Interest in distance graphs sprouted naturally from the chromatic number of the plane problem, which is the problem of determining $\chi(\mathbb R^2, 1)$, in which the distance on $\mathbb R^2$ is the usual Euclidean distance. This problem was formulated in 1950 by Edward Nelson, then an 18-year-old student at the University of Chicago. For the fascinating story of how this problem became known and led to amazing developments in geometric combinatorics, see \cite{s1} and \cite{s2}.

In particular,  coloring distance graphs involving multiple distances is one of the important research problems in the literature of geometric combinatorics. People are still interested in this problem, see for example \cite{cjw} which also contains an extensive literature review for coloring distance graphs involving multiple distances  and see for example \cite{bb} for the same problem for integer distance graphs. In our article we study  coloring integer distance graphs involving multiple distances from a different point of view. So to say, we study a sort of ``worst case scenario" of how you chose your distances. We explain below what we mean.

If $ \left\Vert \cdot  \right\Vert$ is a norm on $\mathbb R^n$, then, clearly, $\chi((\mathbb R^n, \left\Vert \cdot  \right\Vert), d) = \chi((\mathbb R^n, \left\Vert \cdot  \right\Vert), 1)$ for all $d > 0$, so $B_1(\mathbb R^n, \left\Vert \cdot  \right\Vert) = \chi((\mathbb R^n, \left\Vert \cdot  \right\Vert), 1)$. In retrospect, one might think that once we have distance graphs the Babai numbers would be quite natural to think of. Yet, they apparently occured only in the 1990's, to Babai and, possibly, Erd\H{o}s. We conjecture that because of the roles of distance 1 and the Euclidean plane in the problem that started it all, the minds of geombinatorialists were hypnotically distracted for 30 or 40 years.

In this paper we leave the world of real normed spaces and continue the work in \cite{jk} in a quite different class of metric spaces.

Suppose $H$ is a connected simple graph. Let $dist_H$ denote the usual distance function in $H$: for $u, v \in V(H), dist_H(u, v)$ is the length of (the number of edges in) a shortest walk in $H$ from $u$ to $v$, or from $v$ to $u$. (Such a shortest walk will always be a path.) Then $(V(H), dist_H)$ is a metric space, to which we can assign the Babai numbers. The $kth$ Babai number of $(V(H), dist_H)$  will be denoted $B_k(H)$.

If $V(H)$ is finite, then $dist_H$ takes only the values $0,1,\ldots,diam(H)$, where $diam(H)=max[dist_H(u,v) : u,v \in V(H)]$ is the \textit{diameter} of $H$.

Inspection of the definitions reveals that there is no need to confine the discussion of $B_k(H)$ to connected simple graphs $H$. In a disconnected graph $H$, $dist_H(u, v)$ is still defined for vertices $u$ and $v$ in the same component of $H$. The meanings of ``\textit{a distance $d$ is forbidden for a color $c$ in a coloring $\varphi$}'' and ``\textit{a distance $d$ is forbidden by a coloring $\varphi : V(H) \to C$}'' are clear. Under these extended definitions, it is evident that if $H$ has components $H_1, H_2, \ldots, H_t$, then for any positive integer $k$, $B_k(H) = max[B_k(H_j): 1 \leq j \leq t]$.

As announced in the Abstract, our aim here is to evaluate $B_k(P_n)$, where $n \geq 2$, $P_n$ is the path on $n$ vertices, and $1 \leq k \leq n-1 = diam(P_n)$, and also $B_k(C_n)$, where $C_n$ is the cycle on $n$ vertices, and $k \in \{1,2\}$. But we also want to get some answers to the following question: given a connected graph $H$ and $k \in R(V(H), dist_H)=\{1, 2, \ldots, diam(H)\}$, what are the different chromatic numbers $\chi(H, D) = \chi((V(H), dist_H), D)$ as $D$ ranges over $k$-subsets of $\{1, 2, \ldots, diam(H)\}$?

\begin{definition}
Suppose that $H$ is a finite connected graph and $k \in \{1,2, \ldots,$ $diam(H)\}$. The set
$$\{\chi(X, D) \mid D \subseteq \{1,2, \ldots, diam(H)\} \text{ and } \vert D \vert = k\}$$
is the ``\textit{Babai $k$-spectrum}'' of the graph $H$, and is denoted by $Spec(H, k)$.
\end{definition}

From here on, all graphs will be finite and simple.

\section{Babai Numbers and Babai Spectra of Paths}

\begin{theorem}
\label{path}
For all positive integers $n$ and $k$ such that $n \geq k+1$, we have:
$$B_k(P_n)=k+1.$$    
\end{theorem}

\begin{proof} 
We give two proofs. One of them uses the well-known greedy or Grundy algorithm with reference to the natural order of the vertices along the path; and the other uses an algorithmic approach applying the division algorithm. \\
Suppose $n \geq k+1$, and let the vertices of $P_n$ from one end to the other be $v_1,v_2, \ldots, v_n$. Clearly, $B_k(P_n) \geq \chi(P_n, \{ 1,2, \ldots, k\}) \geq k+1$, because the vertices $v_1, v_2, \ldots, v_{k+1}$ induce a complete graph $K_{k+1}$ in $G(P_n, \{ 1,2, \ldots, k\})$. Let $D \subseteq \{1,2, \ldots, diam(P_n)=n-1\}$, $|D|=k$. \\
{\it First Proof.}
We will properly color the graph $G(P_n, D)$ with no more than $k+1$ colors by starting at $v_1$ and coloring the $v_i$ ``greedily'' with positive integers: Color $v_1$ with 1. Having properly colored $v_1, v_2, \ldots, v_{i-1}$, $i \leq n$, with integers, color $v_i$ with the smallest integer not appearing as a color on neighbors of $v_i$ in $G(P_n, D)$ among $v_1, v_2, \ldots, v_{i-1}$. Since $v_i$ can have no more than $k$ neighbors among $v_1, v_2, \ldots, v_{i-1}$, the color on $v_i$ will be $\leq k+1$. So $B_k(P_n)=k+1$. \\
{\it Second Proof.} Let $D=\{d_1, d_2, \ldots, d_k\} \subseteq \{1,2, \ldots, n-1\}$ be a given set of distances such that $d_i > d_{i+1}$ for each $i=1,2, \ldots, k-1$. By division algorithm, for each $i=1,2, \ldots, k-1$, there exist unique $q_i$ and $r_i$ such that $d_i=q_i d_{i+1}+r_i$, where $0 \leq r_i < d_{i+1}$. Again, we will properly color the graph $G(P_n, D)$ with no more than $k+1$ colors by starting at $v_1$ with positive integers: Color the first $d_k$ vertices by the color 1, the second $d_k$ vertices by the color 2, the third $d_k$ vertices by the color 1, and continue coloring alternately with 1 and 2 up to and including the first $q_{k-1}d_k$ vertices, then color the next $r_{k-1}$ vertices by color 1 if the previous $d_k$ vertices are colored with 2, otherwise color with 2. Thus, the first $d_{k-1}$ vertices are colored with at most two colors such that the distance $d_k$ is forbidden. Then, color the second $d_{k-1}$ vertices by the same approach, by using color 3 instead of color 1, color 1 instead of color 2.  Continue coloring alternately up to and including the first $q_{k-2}d_{k-1}$ vertices, then color the next $r_{k-2}$ vertices with the same colors used to color the penultimate $d_{k-1}$ vertices, starting from the first vertex of these $d_{k-1}$ vertices. Thus, the first $d_{k-2}$ vertices are colored with at most three colors such that the distances $d_k$ and $d_{k-1}$ are forbidden. Then, color the second $d_{k-2}$ vertices by the same approach, by using color 4 instead of color 1, color 1 instead of color 2, and color 2 instead of color 3. Continue coloring alternately up to and including the first $q_{k-3}d_{k-2}$ vertices, then color the next $r_{k-3}$ vertices with the same colors used to color the penultimate $d_{k-2}$ vertices, starting from the first vertex of these $d_{k-2}$ vertices. Thus, the first $d_{k-3}$ vertices are colored with at most four colors such that the distances $d_k$, $d_{k-1}$ and $d_{k-2}$ are forbidden. By continuing in this way, after a finite number of steps, all the vertices of $G(P_n, D)$ are properly colored with at most $k+1$ colors in such a way that all the distances $d_1, d_2, \ldots, d_k$ are forbidden.
\end{proof}

\begin{lemma}
Let $n>m>1$ be positive integers, and let $(m,n)$ and $(1,n)$ denote the open intervals.
\begin{enumerate}
\item The number of integers divisible by $m$ in $(m,n)$ is equal to $\left\lceil \frac{n}{m} \right\rceil - 2$.
\item The number of integers divisible by $m$ in $[1,n)$ is equal to $\left\lceil \frac{n}{m} \right\rceil - 1$.
\end{enumerate}    
\end{lemma}

\begin{proof}
Let $q=\left\lfloor \frac{n}{m} \right\rfloor \geq 1$ and $0 \leq r < m$ be such that $n=qm+r$. If $r=0$ then the multiples of $m$ strictly between $m$ and $n$ are $2m, 3m, \ldots, (q-1)m$, so the number of them is $q-2= \left\lfloor \frac{n}{m} \right\rfloor-2 = \left\lceil \frac{n}{m} \right\rceil-2$. If $r > 0$ then the multiples of $m$ strictly between $m$ and $n$ are $2m, 3m, \ldots, qm$, so the number of them is $q-1= \left\lfloor \frac{n}{m} \right\rfloor-1 = (\left\lceil \frac{n}{m} \right\rceil-1)-1 = \left\lceil \frac{n}{m} \right\rceil-2$. Similarly, (2) holds.
\end{proof}

\begin{theorem}
Let $n \geq k+1$ be fixed positive integers such that $k \leq \left\lfloor \frac{n}{2} \right\rfloor$. Then: 
$$Spec(P_n, k)=\{2, 3, \ldots, k+1\}.$$
\end{theorem}

\begin{proof}
First of all, since $k \leq \left\lfloor \frac{n}{2} \right\rfloor$, we can form a $k$-subset $D$ of $\{1,2, \ldots, n-1\}$ with every element odd which gives $\chi(P_n, D)=2$: just color odd numbered vertices with the color 1 and the even numbered vertices with the color 2. Secondly, from the proof of Theorem \ref{path}, for $D=\{1,2, \ldots, k\}$, $\chi(P_n, D)=k+1$.
\hspace*{10pt}Now, let $m \in \{3, 4, \ldots, k\}$. Choose $D=\{1, 2, \ldots, m-1, d_1, d_2, \ldots, d_t\}$, where $t=k-m+1$, and for each $i=1, 2, \ldots, t$, $m< d_i < n$ is not a multiple of $m$. (There are $n-m-1$ positive integers in the open interval $(m,n)$, and $\left\lceil \frac{n}{m} \right\rceil -2$ of them are divisible by $m$. Therefore, the number of positive integers which are not multiples of $m$ in the open interval $(m,n)$ is equal to $(n-m-1)-\left(\left\lceil \frac{n}{m} \right\rceil -2\right)=(n-m+1)-\left\lceil \frac{n}{m} \right\rceil \geq (n-m+1)-\left\lceil \frac{n}{2} \right\rceil = \left\lfloor \frac{n}{2} \right\rfloor - m + 1 \geq k-m+1=t$. Thus, there are at least $t$ positive integers which are not multiples of $m$ in the open interval $(m,n)$.) Then, $\chi(P_n, D)=m$: The distances $1,2, \ldots, m-1$ gives us a $K_m$ in $G(P_n, D)$, and thus $\chi(P_n, D) \geq m$. Now, color the vertices of $P_n$ with $m$ colors with the coloring
$\varphi: V(P_n)=[n] \to \{1,2, \ldots, m\}$ such that  $\varphi(i)=i \pmod{m}$.
Then, the vertices of $G(P_n, D)$ have the same color if and only if the distance between them is a multiple of $m$. Therefore, all the distances $d_1,d_2, \ldots, d_t$ are automatically forbidden with respect to this coloring. And thus, $\chi(P_n, D) \leq m$. As a result, $\chi(P_n, D) = m$. 
\end{proof}

\begin{theorem}
Let $n, k$ be fixed positive integers such that $\left\lfloor \frac{n}{2} \right\rfloor < k \leq n-1$, and $m \in \{2,3, \ldots, n-1 \}$ be the unique integer such that $n-\left\lceil \frac{n}{m} \right\rceil < k \leq n-\left\lceil \frac{n}{m+1} \right\rceil$. Then
$$\{m+1, m+2, \ldots, k+1\} \subseteq Spec(P_n, k).$$
\end{theorem}

\begin{proof}
Since the (sometimes empty) intervals $J(m)=\left( n-\left\lceil \frac{n}{m} \right\rceil , n-\left\lceil \frac{n}{m+1} \right\rceil \right]$, $m = 2,3, \ldots, n-1$, partition the interval $\left(\left\lfloor \frac{n}{2} \right\rfloor , n-1\right]$, it is clear that for each integer $k$ in that interval there is a unique value of $m$ such that $k$ is in $J(m)$.\\
\hspace*{10pt}Now, since the number of integers in the half open interval $[1,n)$ which are divisible by $m+1$ is $\left\lceil \frac{n}{m+1} \right\rceil-1$, those which are not divisible by $m+1$ is $n- \left\lceil \frac{n}{m+1} \right\rceil$. Now, since $k \leq n-\left\lceil \frac{n}{m+1} \right\rceil$, we can choose $k$ elements from the set $S=\{1,2, \ldots, n-1\}$ such that none is divisible by $m+1$. Thus, for $m+1 \leq s \leq k$, since $n-\left\lceil \frac{n}{m+1} \right\rceil \leq n-\left\lceil \frac{n}{s} \right\rceil$, we can again choose $k$ elements from the set $S$, which are not divisible by $s$. Therefore, for each $s \in \{m+1, m+2, \ldots, k\}$, for $D_s=\{ 1,2, \ldots, s-1, d_1, d_2, \ldots, d_t \}$, where $t=k-s+1$, and $d_1, d_2, \ldots, d_t$ are distinct integers in $\{s+1, \ldots, n-1\}$ which are not multiples of $s$, $\chi(P_n, D_s)=s$ as in the proof of the previous theorem. Also, obviously, for $D_{k+1}=\{1, 2, \ldots, k\}$, $\chi(P_n, D_{k+1})=k+1$. As a result, $Spec(P_n, k) \supseteq \{m+1, m+2, \ldots, k+1\}$.
\end{proof}

All the numerical examples we have done give us the courage to make the following conjecture.

\begin{conjecture}
\textit{Under the hypothesis of the previous theorem:}
$$Spec(P_n, k)=\{m+1, m+2, \ldots, k+1\}.$$
\end{conjecture}

\section{Babai Numbers and Babai Spectra of Cycles}

\counterwithin{lemma}{subsection}\setcounter{lemma}{\thelemmaaux}
\counterwithin{theorem}{subsection}\setcounter{theorem}{\thetheoremaux}
\counterwithin{definition}{subsection}\setcounter{definition}{\thedefinitionaux}

\hspace*{10pt}Throughout this section, the vertices of an $n$-cycle $C_n$ will be denoted usually by the integers $0,1,\ldots,n-1$, and an algebraic structure will be put on $V(C_n)$ by using the correspondence between $V(C_n)$ and $\mathbb{Z}_n=\{0,1,\ldots,n-1\}$. Thus, $i$ and $j$ are adjacent on the cycle if and only if $\vert i-j \vert = 1 \pmod{n}$. Also, since the diameter of $C_n$ is $\left\lfloor \frac{n}{2} \right\rfloor$, the set of all possible distances is $S=\{1,2,\ldots, \left\lfloor \frac{n}{2} \right\rfloor\}$.

The crucial part in the correspondence $V(C_n) \leftrightarrow \mathbb{Z}_n$ is the following: For each $d \in S$, the graph $G(C_n, d)$ consists of disjoint union of cycles $C_m$, where $m=\frac{n}{gcd(n, d)}$, the order of $d$ in the additive group $\mathbb{Z}_n$. And each such cycle $C_m$ corresponds to a coset of the subgroup $\langle d \rangle$ generated by $d$ in the group $\mathbb{Z}_n$. Worth nothing: $\langle d \rangle=\langle gcd(n,d) \rangle$. In the graph $G(C_n, d)$, since there is no edge between any two vertices belonging to two different cosets, for a proper coloring of the vertices of this graph, we only need to consider coloring of vertices in one coset, $\langle d \rangle$; any coloring of $\langle d \rangle$ can be copied on each of $\langle d \rangle$'s other cosets. 

The group $\mathbb{Z}_n$ has a unique subgroup of order $d$ for each positive divisor $d$ of $n$, namely, $\langle n/d \rangle$.

\subsection{The First Babai Number and Babai 1-Spectra of $\pmb{C_n}$}
\begin{theorem}
\label{cycle1}
For all positive integers $n \geq 3$, we have:
$$B_1(C_n)= 
\begin{cases}
2, \text{ if } n=2^k \text{ for some $k \in$} \mathbb{Z^+}; \\
3, \text{ otherwise}.
\end{cases}$$
Also, if $n \neq 2^k$ for any $k \in \mathbb Z^+$, then $Spec(C_n, 1)=\{2,3\}$ if $n$ is even and $Spec(C_n, 1)=\{3\}$ if $n$ is odd.
\end{theorem}

\begin{proof}
Let $D=\{d\}$, where $1 \leq d \leq \left\lfloor \frac{n}{2} \right\rfloor$. The vertices of $\langle d \rangle$ induce a cycle of length $\frac{n}{gcd(n,d)}$ in $G(C_n, d)$. The conclusions of the theorem follow from previous remarks and the well-known facts that even cycles have chromatic number 2 and odd cycles have chromatic number 3. 
\end{proof}

Before undertaking the determination of $B_2(C_n)$ and $Spec(C_n,2)$ for all $n \geq 4$, it will be helpful to develop some machinery in a more general setting. 

\subsection{Useful Facts About Cayley Graphs}

\hspace*{10pt}Suppose that $(A, +, 0)$ is an abelian group and $S \subseteq A \setminus \{0\}$. The \textit{Cayley graph} $Cay(A, S)$ is the graph with vertex set $A$ and $x,y \in A$, $x \neq y$, adjacent if and only if $y-x \in S \cup (-S)$. (In some definitions, $S$ is required to generate $(A,+,0)$, that is, $\langle S \rangle = A$. That will not required here.) The components of our $Cay(A, S)$ will be induced by the cosets of $\langle S \rangle$ in $A$. Clearly, $\chi(Cay(A, S)) = \chi(Cay(\langle S \rangle, S))$.

Suppose that $n \geq 4$ is an integer, $0 < k \leq \left\lfloor \frac{n}{2} \right\rfloor$ is an integer, $D \subseteq \{1, \ldots, \left\lfloor \frac{n}{2} \right\rfloor\}$, and $\vert D \vert = k$. Then the distance graph $G(C_n, D)$ is also the Cayley graph $Cay(\mathbb{Z}_n, D)$; in the latter, the elements of $D$ are considered to represent congruence classes of integers modulo $n$, whereas in the former they are ordinary positive integers, possible distances between vertices on the cycle $C_n$. But although the meaning of $D$ changes in the graph names, the graphs are the same.

\begin{definition}
Suppose that $(A, +, 0)$ is an abelian group, $S \subseteq A \setminus \{0\}$, and $r > 1$ is an integer. $S$ is \textit{weakly $r$-free} if and only if $s_1, s_2, \ldots, s_k \in S$, $m_1, m_2, \ldots, m_k \in \mathbb{Z}$, and $\sum \limits_{j=1}^k m_js_j = 0$ imply that $\sum \limits_{j=1}^k m_j \equiv 0 \pmod{r}$.   
\end{definition}

On the naming of this property: ``free" is often an alternative for ``independent" in algebra. If we were to define \textit{strong $r$-freeness} of a set $S \subseteq A \setminus \{0\}$ as linear independence of sets of vectors in vector spaces is defined, we would require that $\sum \limits_{j=1}^k m_js_j = 0$ only when $m_j \equiv 0 \pmod{r}$ for each $j \in \{1,2, \ldots, k\}$. Thence the word ``weakly" in the definition.

\begin{lemma}
\label{Wrf}
Suppose that $(A, +, 0)$, $S \subseteq A \setminus \{0\}$, and $r > 1$ are as in the definition above, and $S$ is weakly $r$-free. Suppose that $m_1, m_2, \ldots, m_k$, $m'_1, m'_2, \ldots, m'_k \in \mathbb{Z}$, $s_1,s_2, \ldots, s_k \in S$ and $\sum \limits_{j=1}^k m_js_j = \sum \limits_{j=1}^k m'_js_j$. Then \linebreak 
$\sum \limits_{j=1}^k m_j \equiv \sum \limits_{j=1}^k m'_j \pmod{r}$.
\end{lemma}

The proof is straightforward.

\begin{theorem}
\label{Cayley}
Suppose that $(A, +, 0)$ is an abelian group, $S \subseteq A \setminus \{0\}$, $r > 1$ is an integer, and $S$ is weakly $r$-free. Then $\chi(Cay(A, S)) \leq r$.
\end{theorem}

\begin{proof}
Since (as alluded earlier) proper coloring of $Cay(\langle S\rangle, S)$ can be copied on each coset of $\langle S\rangle$ in $A$ so that the coset is properly colored in $Cay(A, S)$ (when $A$ is infinite the copying may have to be assisted by the Axiom of Choice, which will allow the choice of coset representatives), it will suffice to show that $\chi(Cay(\langle S\rangle, S)) \leq r$. \\
\hspace*{10pt}Each of $a \in \langle S\rangle$ is representable as an integer combination $a=\sum \limits_{j=1}^k m_js_j$ in which $m_1, m_2, \ldots, m_k \in \mathbb{Z}$, $s_1,s_2, \ldots, s_k \in S$. Color such an $a$ with the congruence class of $\sum \limits_{j=1}^k m_j \pmod{r}$. The coloring is well defined by Lemma \ref{Wrf}. It is proper because of the way adjacency is defined in $Cay(A, S)$. 
\end{proof}

\begin{lemma}
\label{Dwrf}
Suppose that $n>3$ and $r>1$ are integers, $r \mid n$, and $D=\{d_1, d_2, \ldots, d_k\}$, where $1 \leq d_1 < d_2 < \ldots < d_k \leq n-1$. Suppose that, for some $q \in \{1, 2, \ldots, r-1\}$, $gcd(r, q)=1$ and $d_j \equiv q \pmod{r}$, $j=1,2, \ldots, k$. Then $D$ is weakly $r$-free.
\end{lemma}

\begin{proof}
Any integer combination of the $d_j$ can be rewritten without changing its congruency modulo any integer $> 1$ as $\sum \limits_{j=1}^k m_jd_j$, $m_1, m_2, \ldots, m_k \in \mathbb{Z}$. If $\sum \limits_{j=1}^k m_jd_j \equiv 0 \pmod{n}$, then $\sum \limits_{j=1}^k m_jd_j \equiv 0 \pmod{r}$, because $r \mid n$. Therefore, $(\sum \limits_{j=1}^k m_j) q \equiv \sum \limits_{j=1}^k m_j q \equiv \sum \limits_{j=1}^k m_jd_j \equiv 0 \pmod{r}$, because $d_j \equiv q \pmod{r}$ for each $j=1,2, \ldots, k$; and thus $\sum \limits_{j=1}^k m_j \equiv 0 \pmod{r}$, because $gcd(r,q)=1$.  
\end{proof}

\begin{lemma}
\label{LastLemma}
Suppose that $n > r > 1$ are integers and $r \mid n$. Suppose $k > 1$ is an integer satisfying $(k-1)r+1 \leq \left\lfloor \frac{n}{2} \right\rfloor$, $D = \{tr+1 \mid t \in \{0,1, \ldots, k-1 \}\}$. Then $\chi(C_n, D) \leq r$.
\end{lemma}

\begin{proof}
Since $gcd(r,1)=1$, the conclusion follow from Lemma \ref{Dwrf} and Theorem \ref{Cayley}.
\end{proof}

\begin{corollary}
Suppose that $n,k$ are positive integers, $3 \mid n$, $k > 1$, $n$ is odd, and $3(k-1)+1 \leq \left\lfloor \frac{n}{2} \right\rfloor$. Then $3 \in Spec(C_n, k)$.
\end{corollary}

\begin{proof}
Let $D = \{ 3t+1 \mid t \in \{0,1, \ldots, k-1 \}\}$, as in Lemma \ref{LastLemma}. By that Lemma, $\chi(C_n, D) \leq 3$. On the other hand, $1 \in D$ implies that $C_n$ is a subgraph of $G(C_n, D)$, so $\chi(C_n, D) \geq 3$, since $n$ is odd.
\end{proof}

We took the easy way out in concluding $\chi(C_n, D) \geq 3$. In fact, if $n \geq 3$ is odd, $\emptyset \neq D \subseteq \{1, 2, \ldots, \left\lfloor \frac{n}{2} \right\rfloor\}$, then $\chi(C_n, D) \geq 3$ because $G(C_n, D)$ will contain a cycle corresponding to a cyclic subgroup of $\mathbb{Z}_n$, which will have order dividing $n$ and, therefore, odd.

Although the results immediately preceding may have their uses, when it comes to coloring distance graphs $G(C_n, D)=Cay(\mathbb{Z}_n, D)$, the following will definitely be a tool worth saving.

\begin{proposition}
Suppose that $n > r> 1$ are integers, $r \mid n$, and $D \subseteq \{1,2, \ldots, \left\lfloor \frac{n}{2} \right\rfloor \}$. 
Suppose that $d \not\equiv 0 \pmod r$ for each $d \in D$.
Then we have: $\chi(Cay(\mathbb Z_n, D))\leq r.$
\end{proposition}

\begin{proof}
Color each $q \in \{ 0,1, \ldots, n-1\}$ with its congruence class modulo $r$. Because $r \mid n$ this defines a valid coloring of $\mathbb Z_n$, and for each $x\in \mathbb Z_n$, $d \in D$, $x \not\equiv x+d \pmod r$.
\end{proof}

\subsection{The Second Babai Number and Babai 2-Spectra of $\pmb{C_n}$}
\hspace*{10pt}Let $D=\{s,t\}$, where $1 \leq s , t \leq \left\lfloor \frac{n}{2} \right\rfloor$, $s, t, n \in \mathbb{Z}$, $s \neq t$, and $n \geq 4$. (Note: $C_3$ has no $2^{nd}$ Babai number.) Then the graph $G(C_n, D)$ is 4-regular or 3-regular if $n$ is even and $t=\frac{n}{2}$, and thus its maximum degree $\Delta =4$ or $\Delta =3$. Therefore, by Brooks' theorem, $\chi(C_n, D)=5$ (respectively, $\chi(C_n, D)=4$) if $G(C_n, D)$ contains a $K_5$ (respectively, a $K_4$), and $\chi(C_n, D) \leq 5$ (respectively, $\chi(C_n, D) \leq 4$) otherwise. Thus, we obtain the following lemma.

\begin{lemma}
\label{B2Cn}
For any $n \geq 4$, $B_2(C_n) \leq 5$.  
\end{lemma}

\begin{theorem}
\label{K5}
There exist integers $s,t$ satisfying $1 \leq s < t \leq \left\lfloor \frac{n}{2} \right\rfloor$ such that the graph $G(C_n, \{s,t\})$ contains a copy of $K_5$ if and only if 5 divides $n$.  
\end{theorem}

\begin{proof}
Suppose that $G(C_n, D)$ contains a $K_5$ for some $D=\{s,t\}$. Let the vertices of $K_5$ be $v_0, v_1, v_2, v_3, v_4$ in clockwise ordering. Firstly, since $G(C_n, D)$ is a 4-regular graph, the vertices of $K_5$ have no other neighbors in $G(C_n, D)$. Secondly, since there are only two distances involved between the vertices of $C_n$ to determine $G(C_n, D)$, then, letting $d$ denote the distance in $C_n$, $d(v_0, v_1)=d(v_1, v_2)=d(v_2, v_3)=d(v_3, v_4)=d(v_4, v_0)=s$. Therefore, since the distances are measured over $C_n$, the number of vertices of $C_n$, namely, $n=5s$. That is, $n$ is divisible by 5. \\
\hspace*{10pt}Conversely, suppose that $n=5m$ for some $m \in \mathbb{Z}^+$. Let $D=\{m, 2m\}$. Then, $G(C_n, D)$ contains a $K_5$: Just choose, in clockwise order, vertices $v_0, v_1, v_2, v_3, v_4$ of $C_n$ such that $d(v_i, v_{i+1})=m$, $i=0,\ldots, 4$. 
\end{proof}

\begin{corollary}
\label{B2=5}
Suppose that $4 \leq n \in \mathbb{Z}$. Then $B_2(C_n)=5$ if and only if 5 divides $n$.
\end{corollary}

\begin{proof}
$B_2(C_n)=5$ if and only if there is some $D=\{s,t\}$, where $1 \leq s < t \leq \left\lfloor \frac{n}{2} \right\rfloor$, such that $\chi(C_n, D)=5$ since for any such $D$, $\Delta(C_n, D) \leq 4$. By Brooks' theorem $\chi(C_n, D)=5$ holds if and only if some component of $G(C_n, D)$ is $K_5$. By Theorem \ref{K5}, there is such a $D$ if and only if 5 divides $n$.
\end{proof}

\begin{lemma}
\label{4yes5no}
If 4 divides $n$ and 5 does not divide $n$, where $4 \leq n \in \mathbb{Z}$, then $B_2(C_n)=4$.
\end{lemma}

\begin{proof}
Let $n=4m$ for some $m \in \mathbb{Z}^+$. Let $D=\{m,2m\}$, and let the vertices of $C_n$ be labeled $v_0, v_1, \ldots, v_{n-1}$, ordered clockwise from $v_0$. Then $v_0, v_m, v_{2m}, v_{3m}$ induce a $K_4$ in $G(C_n, D)$. Therefore $B_2(C_n) \geq 4$. But, also, $B_2(C_n) < 5$, by Corollary \ref{B2=5}. Therefore, $B_2(C_n)=4$.  
\end{proof}

\begin{lemma}
\label{3and5no}
If $4 \leq n \in \mathbb{Z}$ and neither 3 nor 5 divides $n$, then $B_2(C_n)=4$.
\end{lemma}

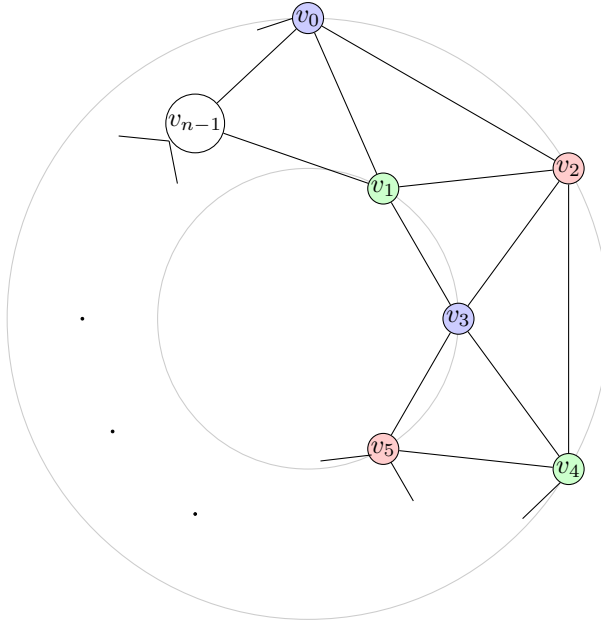
\begin{figure}[htbp]
\centering
\begin{tikzpicture}[inner sep=0.1mm]

\draw (0,0) [draw=black!20, thin] circle [radius=4 cm];
\draw (0,0) [draw=black!20, thin] circle [radius=2 cm];

\foreach \pt/\col/\r/\ang in {$v_0$/blue!20/4/90, $v_1$/green!20/2/60, $v_2$/red!20/4/30, $v_3$/blue!20/2/360, $v_4$/green!20/4/330, 
                              $v_5$/red!20/2/300, $v_{n-1}$/white/3/120} 
                             {\node[circle,draw,fill=\col] (\pt) at (\ang:\r){\pt};}
\foreach \x/\y in {$v_0$/$v_1$, $v_0$/$v_2$, $v_1$/$v_2$, $v_1$/$v_3$, $v_2$/$v_3$, $v_2$/$v_4$, $v_3$/$v_4$, $v_3$/$v_5$, $v_4$/$v_5$,
                   $v_{n-1}$/$v_0$, $v_{n-1}$/$v_1$} 
                  {\draw (\x) -- (\y);}

\filldraw [black] (240:3cm) circle (0.5pt);
\filldraw [black] (210:3cm) circle (0.5pt);
\filldraw [black] (180:3cm) circle (0.5pt);

\draw (93:4) -- (100:3.9);

\draw (327:4) -- (317:3.9);

\draw (300:2.2) -- (300:2.8);
\draw (295:2) -- (275:1.9);

\draw (128:3) -- (134:2.5);
\draw (128:3) -- (136:3.5);

\end{tikzpicture}
\caption{Coloring the vertices of $G(C_n, D_0)$ with three colors is impossible}
\label{fig: 1}
\end{figure}

\begin{proof}
By results above, $5 \nmid n$ implies that $B_2(C_n) \leq 4$. Let $V(C_n)=\{v_0, v_1,$$\ldots,$ $v_{n-1}\}$, ordered clockwise. Let $D_0=\{1, 2\}$. Then, $G(C_n, D_0)$ is a graph consisting of $n$ adjacent triangles sharing an edge. ($G(C_n, D_0)$ can be visualized easily by drawing two concentric circles and locating the vertices having the same index parity on one of these circles and the remaining vertices on the other.)
Now, let us start coloring the vertices of $G(C_n, D_0)$ by first coloring the vertices of the  triangle with vertices $v_0, v_1, v_2$ by the colors 0, 1, 2, say, blue, green, red, respectively (Figure~\ref{fig: 1}). (In Figure~\ref{fig: 1}, $v_{n-1}$ is on the inner (respectively, outer) circle if $n-1$ is odd (respectively, even).) \\
\hspace*{10pt}If $\chi(C_n, D_0)=3$, then this coloring of $v_0, v_1, v_2$ forces a proper coloring $\varphi: \{v_0, v_1, \ldots, v_{n-1}\} \to \{0,1,2\}$. Reading subscripts modulo 3, it is easy to see that $\varphi(v_{i+3}) = \varphi(v_i)$, $i=0,1, \ldots, n-4$. Since $n \not\equiv 0 \pmod{3}$, 
$n-1 \not\equiv 2 \pmod{3}$, and thus $\varphi(v_{n-1})=0$ or 1. But $v_{n-1}$ is adjacent to both $v_0$ and $v_1$. Therefore, $3 < \chi(C_n, D_0) \leq B_2(C_n) < 5$, so $B_2(C_n)=4$.
\end{proof}

\begin{theorem}
\label{B2=4}
Assume that 5 does not divide $n \geq 4, \: n \in \mathbb{Z}$. If $n$ has an odd prime divisor $p \neq 3$, then $B_2(C_n)=4$.
\end{theorem}

\begin{proof}
Assume that there is an odd prime $p \neq 3$ dividing $n$. Let $n=p^a m$, where $a \geq 1$ and $p$ does not divide $m \geq 1$. Now, if 3 does not divide $n$, then $B_2(C_n)=4$ by Lemma \ref{3and5no}. Thus, we assume that 3 divides $n$. Then, 3 divides $m$. And since $p \geq 7$, $2m < \left\lfloor \frac{n}{2} \right\rfloor$. Let $D_0=\{m, 2m\}$. 
The order of the subgroup $\langle m \rangle$ of $\mathbb{Z}_n$ is $\lvert \langle m \rangle \rvert = \frac{p^a m}{gcd(p^a m, m)}=p^a$. \\
\hspace*{10pt}Let the vertices of $C_n$, in order around the cycle, be $v_0, v_1, \ldots, v_{n-1}$, and let $C_{p^a}$ be the cycle with vertices, around the cycle, $v_0, v_m, v_{2m}, \ldots, v_{(p^a-1)m}$. Since $p$ is odd, we have $\chi(C_n, D_0) \geq \chi(C_{p^a})=3$. Since $B_2(C_n) < 5$, by Corollary \ref{B2=5}, to conclude that $B_2(C_n)=4$ it will suffice to show that $\chi(C_n, D_0) > 3$. \\
\hspace*{10pt}We have $\chi(C_n, D_0) \geq \chi(C_{p^a})=3$, from above. Suppose $G(C_n, D_0)$ is properly colored with colors $0, 1, 2$. Reading subscripts modulo $n$, for each $i \in \{0, 1, \ldots, n-1\}$, $v_i, v_{i+m}, v_{i+2m}$ induce a $K_3$ in $G(C_n, D_0)$. Therefore, if we color $v_{im}$ with $i$ for $i \in \{0, 1, 2 \}$, which we can do without loss of generality, the colors $0, 1, 2$ will have to repeat in that order as we color the vertices $v_{jm}$, $j=0, 1, \ldots, p^a-1$. But then, since $v_0$ is colored 0, for the coloring to be proper, $v_{(p^a-3)m}$ must be colored 0, whence $p^a-3 \equiv 0 \pmod 3$, whence $3 \mid p$, contrary to supposition. Thus $\chi(C_n, D_0) > 3$.
\end{proof}

The results preceding determine $B_2(C_n)$ for all $n \geq 4$ except $n=3^u, u \geq 2$ and $n=2 \mathord{\cdot} 3^u, u \geq 1$. Since, for such $n$, $5 \nmid n$ and, clearly, $G(C_n, \{s,t\})$, where $s,t$ are distinct integers such that $1 \leq s, t \leq \frac{n}{2}$, contains an odd cycle, we have $B_2(C_n) \in \{3,4\}$. We shall show that $B_2(C_n)=3$ for such $n$ by showing that for $n=3^u, u \geq 2$, $Spec(C_n, 2)=\{3\}$ and for $n=2 \mathord{\cdot} 3^u, u \geq 1$, $Spec(C_n, 2)=\{2, 3\}$.

First, a general proposition.

\begin{proposition}
Suppose that $n, k, d_1, d_2, \ldots, d_k \in \mathbb Z^+$, $k>1$, and $1\leq d_1 < d_2 < \ldots < d_k \leq \frac{n}{2}$, and that $g$ is a common divisor of $n, d_1, d_2, \ldots, d_k$. Let $D=\{d_1, d_2, \ldots, d_k\}$ and $\frac{1}{g}D=\{\frac{d_1}{g}, \frac{d_2}{g}, \ldots, \frac{d_k}{g}\}$. Then:
$$\chi(C_n, D) = \chi(C_{n/g}, \frac{1}{g}D).$$
\end{proposition}

\begin{proof}
Let $A$ be the subgroup of $\mathbb Z_n$ generated by $d_1, d_2, \ldots, d_k$ and let $B$ be the subgroup of $\mathbb Z_{n/g}$ generated by $\frac{d_1}{g}, \frac{d_2}{g}, \ldots, \frac{d_k}{g}$. Clearly 
$Cay(A, \{d_1, d_2, \ldots, d_k\}) \simeq Cay(B, \{\frac{d_1}{g}, \frac{d_2}{g}, \ldots, \frac{d_k}{g}\})$. Thus
\begin{align*}
\chi(C_n, D) & =\chi(Cay(A, \{d_1, d_2, \ldots, d_k\})) \\
& =\chi(Cay(B, \{\frac{d_1}{g}, \frac{d_2}{g}, \ldots, \frac{d_k}{g}\}))=\chi(C_{n/g}, \frac{1}{g}D).
\end{align*}
\end{proof}

In what follows we shall be properly coloring $G(C_n, \{s,t\})=Cay(\mathbb Z_n, \{s,t\})$, for $n=3^u, u \geq 2$ and $n=2 \mathord{\cdot} 3^u, u \geq 1$, with 3 colors $0,1,2$, standing for the congruence classes in $\mathbb Z$ modulo ${3}$. 

In defining these colorings we will take advantage of the following elementary convenience: If $n$ is a positive integer divisible by 3, then to define a function from $\mathbb Z_n$ to $\mathbb Z_3$ it suffice to define function from $\{0,1, \ldots, n-1\}$ to $\{0,1,2\}$.

\begin{lemma}
Suppose that $u \geq 2, s,t$ are integers, $n=3^u$, $1 \leq s,t  \leq \frac{n}{2}$, $3 \mid s$, and $3 \nmid t$. Then $\chi(C_n, \{s,t\})=3$.
\end{lemma}

\begin{proof}
Obviously, $G(C_n, \{s,t\})=G$ contains an odd cycle, so $\chi(G) \geq 3$. \\
\hspace*{10pt}We shall use $G=Cay(\mathbb Z_n, \{s,t\})$ and properly color $G$ by a coloring $\varphi: \{0,1, \ldots, n-1\}=\mathbb Z_n \to \{0,1,2\}$, as explained above. Let $s=3^a r,$ $1 \leq a \leq u-1$, where $3 \nmid r$. Then $gcd(n,s)=3^a$ so $\langle s \rangle=3^{u-a}$. \\
\hspace*{10pt}Let $I_s=\{\lambda \in \mathbb Z \mid 0 \leq \lambda \leq 3^{u-a}-1\}$ and $J_t=\{\mu \in \mathbb Z \mid 0 \leq \mu \leq 3^{a}-1\}$, and define $\psi: I_s \times J_t \to \mathbb Z_n$ by $\psi(\lambda, \mu) \equiv \lambda s+ \mu t \pmod{n}$. We claim that $\psi$ is one-to-one, and, therefore, onto, since $\vert I_s \times J_t \vert =3^{u-a}3^a=n$. \\
\hspace*{10pt}To see that $\psi$ is one-to-one, suppose that $0 \leq \lambda, \lambda' \leq 3^{u-a}-1$, $0\leq \mu, \mu' \leq 3^a-1$, and
$$\psi(\lambda, \mu) \equiv \lambda s+ \mu t \equiv \lambda' s + \mu' t \equiv \psi(\lambda', \mu') \pmod{n}.$$ 
Then 
$$(\lambda - \lambda')s \equiv (\mu' - \mu) t \pmod n.$$ 
Therefore, because $3^a \mid s, n$, and $3 \nmid t$, it must be that $3^a \mid \mu' - \mu$. Therefore, since $\vert \mu'-\mu \vert <3^a$, $\mu=\mu'$. From
$(\lambda - \lambda')s \equiv 0 \pmod n$ we obtain that $\lambda - \lambda'=3^{u-a}k$ for some integer $k$. Since $\vert \lambda - \lambda' \vert < 3^{u-a}$, it follows that $\lambda=\lambda'$. \\
\hspace*{10pt}We can now 3-color $\mathbb Z_n$ by defining $\varphi: I_s \times J_t \to \{0,1,2\}$ (with $0,1,2$ standing for the congruence classes $3\mathbb Z, 1+3\mathbb Z, 2+3\mathbb Z$ of the integers modulo 3). The formula for $\varphi$ will depend on $r=\frac{s}{3^a}$ and $t$. Recalling that neither is divisible by 3, we have that $rt \equiv 1$ or 2 $\pmod 3$. \\
If $rt \equiv 1 \pmod 3$, we define 
\begin{align*}
\varphi(\lambda, \mu) \equiv \lambda+\mu \pmod 3; 
\end{align*}
if $rt \equiv 2 \pmod 3$, we define 
\begin{align*}
\varphi(\lambda, \mu) \equiv 2\lambda+\mu \pmod 3.
\end{align*}
\hspace*{10pt}For $y \equiv \lambda s+\mu t \pmod n$, $\lambda \in I_s, \mu \in J_t$ we allow the notation $\varphi(y)=\varphi(\lambda, \mu)$. To show that $\varphi$ is a proper coloring of $G$ it suffices to show that for all $y \in \mathbb Z_n$, $\varphi(y+s)\neq \varphi(y)\neq \varphi(y+t)$. \\
\hspace*{10pt}If $0 \leq \lambda \leq 3^{u-a}-2$, $\mu \in J_t$, then $\lambda+1 \leq 3^{u-a}-1$, i.e., $\lambda + 1 \in I_s$, so, if $y=\lambda s + \mu t$, $y+s=(\lambda + 1)s+\mu t$. This implies that
\begin{flalign*}
\varphi(y+s) &  \equiv q(\lambda+1)+\mu \pmod 3, \: q \in \{1,2\} \\
& \equiv q\lambda+\mu+q
\equiv \varphi(y) + q 
\not\equiv \varphi(y) \pmod 3,
\end{flalign*}
whichever $\varphi$ we are applying ($q=1$ or $q=2$). \\
\hspace*{10pt} What if $\lambda=3^{u-a}-1$, $\mu \in J_t$, and $y=\lambda s+\mu t$? Because $3^{u-a}s \equiv 0 \pmod n$, we have, whether $q=1$ or $q=2$,
\begin{align*}
\varphi(y+s) = \varphi(3^{u-a}s+\mu t) \equiv \varphi(0.s+\mu t) \equiv \mu  \not\equiv -q + \mu  \equiv \varphi(y) \pmod 3.
\end{align*}
\hspace*{10pt}Similarly, if $y=\lambda s + \mu t \pmod n$, $\lambda \in I_s$ and $0 \leq \mu \leq 3^a-2$, then $\varphi(y+t) \neq \varphi(y)$, whichever $\varphi$ we are using. To show that $\varphi$ is a proper coloring, it will suffice to show that if $\lambda \in I_s$ and $y=\lambda s + (3^a-1)t$, then
\begin{align*}
\varphi(y+t) \equiv \varphi(\lambda s + 3^at) \not\equiv \varphi(y) \equiv q \lambda + 3^a - 1 \equiv q\lambda + 2 \pmod 3.
\end{align*}
We have
\begin{align*}
y + t \equiv \lambda s + 3^a t \equiv \lambda_0 s + \mu_0 t \pmod n, 
\end{align*}
for some $0 \leq \lambda_0 \leq 3^{u-a}-1, 0 \leq \mu_0 \leq 3^a -1$. Therefore, for some integer $k$,
\begin{align*}
(\lambda - \lambda_0)s = (\lambda - \lambda_0)3^a r = (\mu_0-3^a)t + k \mathord{\cdot} 3^u.  
\end{align*}
This implies that $3^a \mid (\mu_0-3^a)t$. 
Thus, since $3 \nmid t$,  $3^a \mid (\mu_0-3^a)$. So $3^a \mid \mu_0$. 
Therefore, since $0 \leq \mu_0 < 3^a$, $\mu_0 = 0$. \\
We now have
\begin{align*}
(\lambda -\lambda_0)3^ar = -3^a t + k.3^u.
\end{align*}
And so $(\lambda -\lambda_0)r = -t + k \mathord{\cdot} 3^{u-a}$.
Thus, since $u > a$, $(\lambda -\lambda_0)r \equiv -t \pmod 3$. \\
Since $3 \nmid r$, $r^2 \equiv 1 \pmod 3$; consequently,
\begin{align*}
\lambda - \lambda_0 \equiv -rt \pmod 3 \Rightarrow \lambda \equiv \lambda_0 -rt \pmod 3. 
\end{align*}
We have $\varphi(y) \equiv q \lambda + 2 \equiv q \lambda_0 + 2 -qrt \pmod 3$ while $\varphi(y+t) \equiv q \lambda_0 + \mu_0 \equiv q \lambda_0 \pmod 3$. 
Recollecting that $rt \equiv 1 \pmod 3$ implies that $q=1$ and $rt \equiv 2 \pmod 3$ implies that $q=2$, we see that in either case,
\begin{align*}
\varphi(y + t) \equiv q \lambda_0 \not\equiv q \lambda_0 + 1 \equiv \varphi(y) \pmod 3.
\end{align*}
As a result, $\varphi$ is a proper 3-coloring of $\mathbb Z_n = \{0,1, \ldots, n-1\}$.
\end{proof}

\begin{theorem}
Suppose that $u>1$, $s$ and $t$ are integers, $n=3^u$, $1 \leq s,t \leq \frac{n}{2}$, and $s \neq t$. Then $\chi(C_n, \{s,t\})=3$.
\end{theorem}

\begin{proof}
In all cases, $G=G(C_n, \{s,t\})=Cay(\mathbb Z_n, \{s,t\})$ has odd cycles, so $\chi(G) \geq 3$. \\
{\bf Case 1.} Neither $s$ nor $t$ is divisible by 3. \\
Color $\mathbb Z_n=\{0,1, \ldots, n-1\}$ by $\varphi(y) \equiv y \pmod 3$. The definition is valid because $3 \mid n$; for the same reason, because $s \not\equiv 0 \not\equiv t \pmod 3$, adding in $\mathbb Z_n$ either $s$ or $t$ to any $y \in \mathbb Z_n$ will give an element of $\mathbb Z_n$ with a different color than $y$, so the coloring is proper. \\ 
{\bf Case 2.} Both $s$ and $t$ are divisible by 3. \\
Then $g=gcd(s,t,n)=3^a$ for some $a$, $1 \leq a \leq u-2$; $a \leq u-2$ because $\frac{s}{3^a}, \frac{t}{3^a}$ are distinct positive integers $\leq \frac{3^{u-a}}{2} \cdot$ \\
By Proposition 3.3.1,
\begin{align*}
\chi(C_n, \{s,t\}) = \chi(C_{3^{u-a}}, \{\frac{s}{3^a}, \frac{t}{3^a}\}).    
\end{align*}
If neither $\frac{s}{3^a}, \frac{t}{3^a}$ is divisible by 3, then
\begin{align*}
\chi(C_n, \{s,t\}) = \chi(C_{3^{u-a}}, \{\frac{s}{3^a}, \frac{t}{3^a}\}) \leq 3  
\end{align*}
by Case 1. Otherwise, exactly one of $\frac{s}{3^a}, \frac{t}{3^a}$ is divisible by 3, and $\chi(G) \leq 3$ will follow from Case 3. \\
{\bf Case 3.} Without loss of generality, $3 \mid s$ and $3 \nmid t$. Then $\chi(C_n, \{s,t\})=3$ by Lemma 3.3.4.
\end{proof}

Now, on to $n=2 \mathord{\cdot} 3^u$!

\begin{proposition}
Suppose that $n,k, d_1, d_2, \ldots, d_k \in \mathbb Z^+$, $k > 1$, and $1 \leq d_1 < d_2 < \ldots < d_k \leq \frac{n}{2}$. Suppose that $n$ is even and each $d_i$ is odd. Let $D=\{d_1, d_2, \ldots, d_k\}$. Then $\chi(C_n, D)=2$.
\end{proposition}

\begin{proof}
Color  $y \in \{0,1, \ldots, n-1\}$ with red if $y$ is even and with blue if $y$ is odd. This is a valid coloring of $\mathbb Z_n$ because $n$ is even, and for each $i \in \{1,2, \ldots, k\}$, $y \in \mathbb Z_n$, $y$ and $y+d_i$ will bear different colors because $d_i$ is odd. Therefore, the coloring is a proper coloring of $Cay(\mathbb Z_n, D)$ with 2 colors.
\end{proof}

Proposition 3.3.2 is a special case of the application of Lemma 3.2.2, with $r=2$ and $q=1$, to Theorem 3.2.1 with $A=\mathbb Z_n$, and also of Proposition 3.2.1, with $r=2$.

On to $n=2 \mathord{\cdot} 3^u$!

\begin{lemma}
Suppose that $u \geq 1, s, t$ are integers, $n=2 \mathord{\cdot} 3^u$, $1 \leq s, t\leq \frac{n}{2}$, $6 \mid s$, and neither 2 nor 3 divide $t$. Then 
$\chi(C_n, \{s,t\})=3$.
\end{lemma}

\begin{proof}
For some $a,b$, $1 \leq a \leq u-1$, $1\leq b$, and $r$, which is divisible by neither 2 nor 3, $s=2^b 3^ar$. Then $gcd(s,n)=2 \mathord{\cdot} 3^a$ so if $\langle s \rangle$ is the subgroup generated by $s$ in $\mathbb Z_n$, then $\vert \langle s \rangle \vert = 3^{u-a}$. Therefore, 
$G=G(C_n, \{s,t\})=Cay(\mathbb Z_n, \{s,t\})$ contains an odd cycle, so $\chi(G) \geq 3$. \\
\hspace*{10pt}Let $I_s=\{\lambda \in \mathbb Z \mid 0 \leq \lambda \leq 3^{u-a}-1\}$ and $J_t=\{\mu \in \mathbb Z \mid 0 \leq \mu \leq 2 \mathord{\cdot} 3^{a}-1\}$, define $\psi: I_s \times J_t \to \mathbb Z_n$ by $\psi(\lambda, \mu) \equiv \lambda s+ \mu t \pmod{n}$. As in Lemma 3.3.4, we verify that $\psi$ is one-to-one, and, therefore, onto: Suppose
\begin{align*}
0 \leq \lambda, \lambda' \leq 3^{u-a}-1, 0\leq \mu, \mu' \leq 2 \mathord{\cdot} 3^a-1, \textrm{ and } \lambda s+ \mu t \equiv \lambda' s + \mu' t \pmod{n}.   
\end{align*}
Then 
$$(\lambda - \lambda')s \equiv (\mu' - \mu) t \pmod n.$$ 
Since $2 \mathord{\cdot} 3^a$ divides both $s$ and $n$, and $t$ is divisible by neither 2 nor 3, we have $2 \mathord{\cdot} 3^a \mid \mu' - \mu$. But $\vert \mu'-\mu \vert <2 \mathord{\cdot} 3^a$, so $\mu=\mu'$; thence, from $(\lambda - \lambda')s \equiv 0 \pmod n$, as in the proof of Lemma 3.3.4, we conclude that $\lambda=\lambda'$. \\
\hspace*{10pt}Using the fact for each $y \in \mathbb Z_n$ there is a unique pair $(\lambda, \mu) \in I_s \times J_t$ such that $y \equiv \lambda s+\mu t \pmod n$, we will, as in the proof of Lemma 3.3.4, define $\varphi: \mathbb Z_n \to \{0,1,2\}=\mathbb Z_3$ by
$\varphi(y)=\varphi(\lambda s+\mu t)=q\lambda + \mu \pmod 3$, with $q \in \{1,2\}$ chosen so that $\varphi$, so defined, satisfies, for $y \in \mathbb Z_n$ and $x \in \{s,t\}$, $\varphi(y) \neq \varphi(y+x)$. The value of $q$ will depend on $a,b,r=\frac{s}{2^b 3^a}$, and $t$. \\
\hspace*{10pt}It is straightforward to see that if $0 \leq \lambda \leq 3^{u-a}-2$ and $\mu \in J_t$ then  
\begin{flalign*}
\varphi((\lambda s + \mu t) + s) = \varphi((\lambda + 1) s + \mu) & \equiv q(\lambda + 1) + \mu \pmod 3 \\
& \not\equiv q\lambda + \mu \equiv \varphi(\lambda s + \mu t) \pmod 3.
\end{flalign*}
\hspace*{10pt}Similarly, if $\lambda \in I_s$, $0 \leq \mu \leq 2 \mathord{\cdot}  3^a-2$, and $y=\lambda s + \mu t$, then
\begin{align*}
\varphi(y+t) \equiv q \lambda + \mu + 1 \not\equiv q \lambda + \mu  \equiv \varphi(y) \pmod 3.
\end{align*}
\hspace*{10pt}Therefore, to show that $\varphi$ is a proper coloring of $Cay(\mathbb Z_n, \{s,t\})=G(C_n, \{s,t\})$, it will suffice to show the following: \\
(1) If $\mu \in J_t$ and $y=(3^{u-a}-1)s + \mu t \pmod n$ then $\varphi(y) \neq \varphi(y+s)$; and \\
(2) If $\lambda \in I_s$ and $y=\lambda s + (2 \mathord{\cdot} 3^{a}-1)t \pmod n$ then $\varphi(y) \neq \varphi(y+t)$.

\noindent (1): For some $\lambda_0 \in I_s, \mu_0 \in J_t$,
\begin{align*}
y+s \equiv 3^{u-a}s + \mu t \equiv \lambda_0 s +  \mu_ 0 t \pmod n. 
\end{align*}
This implies that for some integer $k, 3^{u-a}s + \mu t = \lambda_0 s + \mu_ 0 t + 2 \mathord{\cdot} 3^u k$.
Therefore, $(3^{u-a}-\lambda_0)s = (\mu_0 - \mu)t + 2 \mathord{\cdot} 3^u k$.
Since $u > a$, $2 \mathord{\cdot} 3^a \mid s$, and $t$ is divisible by neither 2 nor 3, it follows that $2 \mathord{\cdot} 3^a \mid (\mu_0 - \mu)$. Since $\mu_0, \mu \in J_t$, 
$\vert \mu_0 - \mu \vert < 2 \mathord{\cdot} 3^a$. Therefore, $\mu_0 = \mu$ and we have 
\begin{align*}
(3^{u-a} - \lambda_0)2^b  3^a  r = 2 \mathord{\cdot} 3^u k. 
\end{align*}
And thus $(3^{u-a} - \lambda_0)2^{b-1} r = 3^{u- a}k \equiv 0 \pmod 3$.
Since $2^{b-1} r \not\equiv 0 \pmod 3$, we conclude that $\lambda_0 \equiv 0 \pmod 3$. \\
Consequently,
\begin{align*}
\varphi(y+s) \equiv q\lambda_0 + \mu_0 \equiv 0 + \mu \equiv \mu \pmod 3,
\end{align*}
while
\begin{align*}
\varphi(y) = \varphi((3^{u-a}-1)s + \mu t) \equiv q(3^{u-a}-1) + \mu \equiv 2q + \mu \not\equiv \mu \pmod 3.
\end{align*}

\noindent (2): For some $\lambda_0 \in I_s, \mu_0 \in J_t$,
\begin{align*}
y+t \equiv \lambda s + 2 \mathord{\cdot} 3^{a}t \equiv \lambda_0 s +  \mu_ 0 t \pmod n.
\end{align*}
This implies that for some integer $k, \lambda s + 2 \mathord{\cdot} 3^a t = \lambda_0 s + \mu_ 0 t + 2 \mathord{\cdot} 3^u k$.
Therefore, $(\lambda -\lambda_0)s = (\mu_0 - 2 \mathord{\cdot} 3^a)t + 2 \mathord{\cdot} 3^u k$. Thus $2 \mathord{\cdot} 3^a \mid  \mu_0 - 2 \mathord{\cdot} 3^a$. As a result, $2 \mathord{\cdot} 3^a \mid \mu_0$.
Since $0 \leq \mu_0 < 2 \mathord{\cdot} 3^a$, $\mu_0 = 0$, so we have  
\begin{align*}
(\lambda - \lambda_0)s = -2 \mathord{\cdot} 3^a t + 2 \mathord{\cdot} 3^u k. 
\end{align*}
Thus $(\lambda - \lambda_0)2^{b-1} r = -t +  3^{u- a}k$.
Therefore, since $(2^{b-1} r)^2 \equiv 1 \pmod 3$, $\lambda \equiv \lambda_0 - 2^{b-1} r t \pmod 3$. \\
Then we have
\begin{align*}
\varphi(y) = \varphi(\lambda s + (2 \mathord{\cdot} 3^a-1)t) \equiv q\lambda + 2 \mathord{\cdot} 3^a - 1 & \equiv q (\lambda_0 - 2^{b-1}r t) + 2 \pmod 3 \\
& \equiv q \lambda_0 + (2 - q \mathord{\cdot} 2^{b-1}r t) \pmod 3,
\end{align*}
while
\begin{align*}
\varphi(y+t) \equiv \lambda_0 s + \mu_0 t \equiv q\lambda_0 + 0 \equiv q \lambda_0 \pmod 3.
\end{align*}
Thus $\varphi(y+t) \not\equiv \varphi(y)$ if $2 - q \mathord{\cdot} 2^{b-1}r t  \not\equiv 0 \pmod 3$. \\
If $2^{b-1}  r  t \equiv 1 \pmod 3$, take $q=1$. If $2^{b-1}  r  t \equiv 2 \pmod 3$, take $q=2$.
\end{proof}

\begin{lemma}
Suppose that $u \geq 1, s, t$ are integers, $n=2 \mathord{\cdot} 3^u$, $1 \leq s, t \leq \frac{n}{2}$, $3 \mid s$, $2 \nmid s$, $3 \nmid t$, $2 \mid t$. Then $\chi(C_n, \{s,t\})=3$.
\end{lemma}

\begin{proof}
Since $gcd(t,n)=2$, $\vert \langle t \rangle \vert=3^u$, odd, so $\chi(C_n, \{s,t\}) \geq 3$. \\
\hspace*{10pt}For some $1 \leq a \leq u$, and $r$ relatively prime to $n$, $s=3^ar$. \\
\hspace*{10pt}Let $I_s=\{\lambda \in \mathbb Z \mid 0 \leq \lambda \leq 2 \mathord{\cdot} 3^{u-a}-1\}$ and $J_t=\{\mu \in \mathbb Z \mid 0 \leq \mu \leq 3^{a}-1\}$. 
As in the previous proofs, define $\psi: I_s \times J_t \to \mathbb Z_n$ by $\psi(\lambda, \mu) \equiv \lambda s+ \mu t \pmod{n}$. \\
We verify that $\psi$ is one-to-one, and, therefore, onto, since $\vert I_s \times J_t \vert = n$: \\
Suppose  
\begin{align*}
\lambda, \lambda' \in I_s, \mu, \mu' \in J_t, \textrm{ and }
\lambda s+ \mu t \equiv \lambda' s + \mu' t \pmod{n}.   
\end{align*}
Then 
$$(\lambda - \lambda')s \equiv (\mu' - \mu) t \pmod n ;$$ 
since $3 \nmid t$ and $3^a \mid s,n$, it must be that $3^a \mid \mu'-\mu$. Since $\vert \mu'-\mu \vert <3^a$, it must be that $\mu=\mu'$. From 
$(\lambda - \lambda')s \equiv 0 \pmod n$ and $gcd(r,n)=1$ we obtain $(\lambda - \lambda') 3^a \equiv 0 \pmod n$. It follows that for some integer $k$, 
$\lambda - \lambda'= 2 \mathord{\cdot} 3^{u-a}k$. Then $\vert \lambda - \lambda' \vert < 2 \mathord{\cdot}  3^{u-a}$ implies that $\lambda = \lambda'$.

Again, we can define $\varphi: \mathbb Z_n \to \mathbb Z_3$ by $\varphi(\lambda s + \mu t)=q \lambda + \mu \pmod 3$, for some $q \in \{1,2\}$. Again, $\varphi$ will be a proper coloring of $Cay(\mathbb Z_n, \{s,t\})$ if $q$ is chosen so that 
$\varphi(\lambda s + 3^at) \not\equiv \varphi(\lambda s + (3^a-1)t) \equiv q\lambda + 2  \pmod 3$, if $\lambda \in I_s$, and 
$\varphi(2 \mathord{\cdot} 3^{u-a} s + \mu t) \neq \varphi((2 \mathord{\cdot} 3^{u-a}-1) s + \mu t)$ if $\mu \in J_t$. 
We leave to the reader the pleasure of discovering that $\varphi$ will be proper if we take $q \equiv tr \pmod 3$.
\end{proof}

\begin{theorem}
Suppose that $u, s$, and $t$ are positive integers, $n=2 \mathord{\cdot} 3^u$, $1 \leq s, t \leq 3^u$, and $s \neq t$. Then 
$\chi(C_n, \{s,t\}) \in \{2,3\}$, and both values are possible. 
\end{theorem}

\begin{proof}
If $s$ and $t$ are both odd, $\chi(C_n, \{s,t\}) = 2$ by Proposition 3.3.2. If $s$ and $t$ are both even, then 2 is a common divisor of $s, t$, and $n$. On the one hand, $\frac{n}{gcd(s,n)}$ is odd, so  $G=G(C_n, \{s,t\})$ contains an odd cycle, whence $\chi(G) \geq 3$. On the other hand, by Proposition 3.3.1 and Theorem 3.3.3, $\chi(G) = \chi(C_{3^u}, \{\frac{s}{2}, \frac{t}{2}\})=3$. \\
\hspace*{10pt}Therefore, we may assume that one of $s,t$ is odd and the other is even. Then $\chi(G) \geq 3$. If both are divisible by 3, then $gcd(s,t,n)=3^a$ for some $1 \leq a \leq u-1$. By Proposition 3.3.1, we would then have 
$\chi(G) = \chi(C_{2 \mathord{\cdot} 3^{u-a}}, \{\frac{s}{3^a}, \frac{t}{3^a}\})$. \\
Either neither of $\frac{s}{3^a}, \frac{t}{3^a}$ is divisible by 3 or exactly one is. If exactly one, either it is the one that is even, in which case Lemma 3.3.5 implies what we want, or it is the one that is odd, in which case Lemma 3.3.6 finishes the proof. In case neither of $\frac{s}{3^a}, \frac{t}{3^a}$ is divisible by 3, we have $\chi(G) \leq 3$ by applying Proposition 3.2.1 with $r=3$ to $Cay(C_{2 \mathord{\cdot} 3^{u-a}}, \{\frac{s}{3^a}, \frac{t}{3^a}\})$, and thus conclude that $\chi(G)=3$.
\end{proof}

\begin{theorem}
\label{mainthm}
Let $4 \leq n \in \mathbb{Z}$. Then, the second Babai number of $C_n$ is:
$$B_2(C_n)=
\begin{cases}
3, \text{ if } n=3^u \text{ or } n=2 \mathord{\cdot} 3^u \text{ for some } u \in \mathbb{Z}^+ ; \\
4, \text{ if } n \neq 3^u, n \neq 2 \mathord{\cdot}  3^u, \text{ and } 5 \nmid n ; \\
5, \text{ if } n=5m \text{ for some } m \in \mathbb{Z}^+.
\end{cases}$$
\end{theorem}

\begin{proof} $ $ \newline
\textbf{Case 1.} If $n=3^u$, respectively, $n=2 \mathord{\cdot} 3^u$, for some $u \in \mathbb{Z}^+$, then $B_2(C_n)=3$  by Theorem 3.3.3, respectively, by Theorem 3.3.4. 
\\
\textbf{Case 2.} If $n \neq 3^u, n \neq 2 \mathord{\cdot}  3^u, \text{ and } n \neq 5u \text{ for any } u \in \mathbb{Z}^+$, then $B_2(C_n)=4$ by Theorem \ref{B2=4}. \\
\textbf{Case 3.} If $n=5m \text{ for some } m \in \mathbb{Z}^+$, then $B_2(C_n)=5$ by Corollary \ref{B2=5}.
\end{proof}

\begin{lemma}
If $n>3$ then
$$\chi(C_n, \{1,2\})=
\begin{cases}
3, & \text{if } \: 3 \mid n; \\
4, & \text{if } \: 3 \nmid n \text{ and } n \neq 5.
\end{cases}$$
\end{lemma}

\begin{proof}
Since $G(C_n, \{1,2\})$ contains triangles, $\chi(C_n, \{1,2\}) \geq 3$ for all $n > 3$. \\
\hspace*{10pt}Trying to properly color $G(C_n, \{1,2\})$ with three colors, $r,b,g$, by catch-as-catch-can, going around the cycle coloring the vertices $v_0, v_1, \ldots, v_{n-1}$, we find that if we start by coloring $v_0$ with $r$ and $v_1$ with $b$ then we are forced into a repeating color-word $rbgrbg, \ldots$. If $3 \mid n$ we wind up with color-word $(rbg)^{n/3}$ and a proper 3-coloring of $G(C_n, \{1,2\})$. \\ 
\hspace*{10pt}When $n \equiv 1 \pmod 3$ we are forced, if confined to 3 colors, to the improper coloring $(rbg)^{\frac{n-1}{3}}r$, so that $v_{n-1}$ and $v_0$, at distance 1, bear the same color, r. But if $v_{n-1}$ is colored with a 4th color, $y$, the coloring is proper, so $\chi(C_n, \{1,2\})=4$ when $n \equiv 1 \pmod 3$. Similarly, when $n \equiv 2 \pmod 3$, $\chi(C_n, \{1,2\}) > 3$, and when $n >5$ the word $(rbg)^{\frac{n-5}{3}}rybgy$ shows that $\chi(C_n, \{1,2\})=4$.
\end{proof}

\begin{lemma}
If $n\geq 20$ then $\chi(C_n, \{2,3\})=3$.
\end{lemma}

\begin{proof}
Let the vertices of $C_n$ be $v_0, v_1, \ldots, v_{n-1}$. If $n \geq 7$ then $G(C_n, \{2,3\})$ contains the 5-cycle $v_0v_3v_6v_4v_2v_0$, and so $\chi(C_n, \{2,3\}) \geq 3$. \\
\hspace*{10pt}Consider 3-colorings of $C_n$ obtainable by concatenating the words $rrbbg$ and $rrbbgg$ over the color set $\{r,b,g\}$, and writing the concatenation around the cycle, starting at $v_0$. If $n=5 \lambda + 6 \mu$ for some non-negative integers $\lambda, \mu$ then any concatenation of $\lambda$ of the first word and $\mu$ of the second, in any order, will properly color $G(C_n, \{2,3\})$. By a famous theorem of Frobenius (see \cite{brauer}), if $n \geq (5-1)(6-1) = 20$ then $n$ is so expressible.
\end{proof}

\begin{lemma}
If $n = 2^a q$, where $a \geq 1, q \geq 5, 3 \nmid q$, and $q$ is odd, then $\chi(C_n, \{2,3\})=3$.
\end{lemma}

\begin{proof}
If $n \geq 20$ then the conclusion follows from the previous lemma. Therefore, it suffices to verify the claim for $n=10=2 \mathord{\cdot} 5$ and $n=14=2 \mathord{\cdot} 7$. First of all, since $\vert \langle 2 \rangle \vert = 5$ and $\vert \langle 2 \rangle \vert = 7$ in $\mathbb Z_{10}$ and $\mathbb Z_{14}$ respectively, $G(C_n, \{2,3\})$ contains an odd cycle for both $n=10$ and $n=14$; and thus $\chi(C_n, \{2,3\}) \geq 3$ for both $n=10$ and $n=14$. In fact, 3 colors are enough for proper coloring:  Color the vertices $v_0, v_1, \ldots, v_9$ of $C_{10}$ as follows: $r,r,b,b,g,r,r,b,g,g$; and color the vertices $v_0, v_1, \ldots, v_{13}$ of $C_{14}$ as follows: $r,r,b,b,g,r,r,b,g,g,r,b,b,g$.
\end{proof}

\begin{theorem}
Let $4 \leq n \in \mathbb{Z}$. Then
$$Spec(C_n, 2)=
\begin{cases}
\{3\}, & \text{if } n=3^u \text{ for some } u \in \mathbb{Z}^+; \\
\{4\}, & \text{if } n=4 \text{ or } n=7; \\
\{5\}, &\text{if } n=5; \\
\{2, 3\}, & \text{if } n=2 \mathord{\cdot} 3^u \text{ for some } u \in \mathbb{Z}^+; \\
\{3, 4\}, & \text{if } n=pm, 1 \leq m \in \mathbb{Z} \text{ is odd }, 5 \nmid m, \\
& 5<p \text{ is a prime}, \text{and } n > 7; \\
\{2,3,4\}, & \text{if } n=2^a 3^u \text{ for some } 2 \leq a, 0 \leq u \in \mathbb{Z} \text{ and } n \neq 4; \\
& \text{or } n=pm, 1 < m \in 2\mathbb{Z}, 5 \nmid m, 5 < p \text{ is a prime}; \\
\{3,4,5\}, &\text{if } n=5m \text{ for some odd integer } m \geq 3; \\
\{2,3,4,5\}, &\text{if } n=5m \text{ for some even integer } m \geq 2.
\end{cases}$$
\end{theorem}

\begin{proof} $ $ \newline
\textbf{Case 1.} Assume that $n=3^u \text{ for some } u \in \mathbb{Z}^+$. Then $Spec(C_n,2)=\{3\}$ by Theorem 3.3.3. \\
\textbf{Case 2.} If $n=4$, then there is only one possibility for $D$, namely, $D=\{1,2\}$. It is obvious that $\chi(C_4, D)=4$. If $n=7$, then there are only three possibilities for $D$, namely, $D=\{1,2\}$ or $D=\{1,3\}$ or $D=\{2,3\}$. It can be shown easily that $\chi(C_7, D)=4$ in all these three cases. \\
\textbf{Case 3.} Assume that $n=5$. Then, there is only one possibility for $D$, namely, $D=\{1,2\}$. It is obvious that $\chi(C_5, D)=5$. \\
\textbf{Case 4.} Assume that $n=2 \mathord{\cdot} 3^u \text{ for some } u \in \mathbb{Z}^+$. Then $Spec(C_n,2)=\{2, 3\}$ by Theorem 3.3.4. \\
\textbf{Case 5.} Assume that $n=pm$, where $5<p$ is a prime number and $1 \leq m$ is an odd integer not divisible by 5. First of all, since $B_2(C_n)=4$ by Theorem \ref{mainthm}, $\chi(C_n, D)=4$ for some $D$ satisfying the conditions of the theorem. Secondly, the order of any subgroup of $\mathbb Z_n$ is odd, because $n$ is odd. That is, every coset of every subgroup of $\mathbb Z_n$ corresponds to an odd cycle. Therefore, $\chi(C_n, D) \geq 3$ for any $D$ as described in the statement of the theorem. 
Now, let $n= 2k+1$, where $k \geq 4$, and $D=\{k-1, k\}$. Then $\chi(C_n, D)=3$: 
Let the vertices be $v_0, v_1, \ldots, v_{2k}$ and let the colors be $0,1$, and 2. Color the first $k-1$ vertices with color 0, the following $k-1$ vertices with color 1, and the remaining three vertices with color 2.
Then two vertices have the same color only if the distance between them is strictly less than $k-1$. Thus, there are not any same color vertices at distance $k-1$ or $k$. \\
\textbf{Case 6.} Assume that 5 does not divide $n \neq 4$, and $n=2^a 3^u$, where $a \geq 2$ and $u \geq 0$, or there is a prime $p>5$ such that $n=pm$, where $1 < m \in 2\mathbb{Z}$ and $5 \nmid m$. First of all, since $B_2(C_n)=4$ by Theorem \ref{mainthm}, $Spec(C_n, 2) \subseteq \{2,3,4\}$ and $\chi(C_n, D)=4$ for some $D$ satisfying the conditions of the theorem. \\
\textbf{Subcase a.} Assume that $n=2^a$, where $a \geq 3$.
First, since $n$ is even, $\chi(C_n, \{1,3\})=2$ by Proposition 3.3.2. 
Second, $\chi(C_n, \{1, 2^{a-1}\})=3$: Let the vertices be $v_0, v_1, \ldots$, $v_{2^a-1}$. First of all, for proper coloring, at least 3 colors are needed, because with only 2 colors, $v_0$ and $v_{2^{a-1}}$ must have the same color. In fact, 3 colors are enough for proper coloring: Color $v_j$ with $j$'s congruence class modulo 2,  $0 \leq j \leq 2^{a-1}-1$, color $v_{2^{a-1}}$ with color 2, and then for $2^{a-1} < j < 2^a-1$, color $v_j$ with the congruence class modulo 2 of $j+1$, and finally color $v_{2^a-1}$ with color 2. \\
\textbf{Subcase b.} Assume that $n=2^a 3^u$, where $a \geq 2$ and $u \geq 1$. First, since $n$ is even, $\chi(C_n, \{1,3\})=2$ by Proposition 3.3.2. Second, since $3 \mid n$, $\chi(C_n, \{1,2\})=3$ by Lemma 3.3.7. \\
\textbf{Subcase c.} Assume that $n=pm, 1 < m \in 2\mathbb Z, 5 \nmid m$, and $5 < p$ is prime. First, since $n$ is even, $\chi(C_n, \{1,3\})=2$ by Proposition 3.3.2. Second, if $3 \nmid m$, write $n=2^a q$, where $a \geq 1$, $q > 5$, $5 \nmid q$, $3 \nmid q$, and $q$ is odd, then $\chi(C_n, \{2,3\}) = 3$ by Lemma 3.3.9; and if $3 \mid m$, then $\chi(C_n, \{1,2\})=3$ by Lemma 3.3.7. \\
\textbf{Case 7.} Assume that $n=5m$, where $m \geq 3$ is an odd integer.  
First of all, since $B_2(C_n)=5$ by Theorem \ref{mainthm}, $\chi(C_n, D)=5$ for some $D$ satisfying the conditions of the theorem; and as in Case 5, $\chi(C_n, D) \geq 3$ for any such $D$. Also, again as in Case 5, let $n=2k+1$, where $k \geq 4$, then for $D=\{k-1, k\}$, $\chi(C_n, D)=3$. 
Now, if $3 \nmid m$, then $\chi(C_n, \{1,2\})=4$ by Lemma 3.3.7; and if $3 \mid m$, write $n=3^a z$, where $z \geq 5$ and $3 \nmid z$, then by Proposition 3.3.1 and Lemma 3.3.7,
$\chi(C_n, \{3^a, 2 \mathord{\cdot} 3^a\}) = \chi(C_{n/3^a}, \{1, 2\}) =4$. \\
\textbf{Case 8.} Assume that $n=5m$, where $m \geq 2$ is an even integer. First of all, since $B_2(C_n)=5$ by Theorem \ref{mainthm}, $Spec(C_n, 2) \subseteq \{2,3,4,5\}$ and $\chi(C_n, D)=5$ for some $D$ satisfying the conditions of the theorem. 
First, since $n$ is even, $\chi(C_n, \{1,3\})=2$ by Proposition 3.3.2. 
Second, if  $3 \nmid m$, write $n=2^a q$, where $a \geq 1, q \geq 5, 3 \nmid q$, and $q$ is odd, then $\chi(C_n, \{2,3\})=3$ by Lemma 3.3.9;
and if $3 \mid m$, then $\chi(C_n, \{1,2\})=3$ by Lemma 3.3.7. 
Third, if $3 \nmid m$, then $\chi(C_n, \{1,2\})=4$ by Lemma 3.3.7; and if $3 \mid m$, as in Case 7, write $n=3^a z$, where $z \geq 5$ and $3 \nmid z$, then by Proposition 3.3.1 and Lemma 3.3.7, $\chi(C_n, \{3^a, 2 \mathord{\cdot} 3^a\}) = \chi(C_{n/3^a}, \{1, 2\}) =4$.
\end{proof}

\end{document}